%% file: andkpaper.tex
\documentclass [leqno,notitlepage,12pt]{article}
\usepackage{amssymb,amsmath,amscd}
\usepackage{ifthen}
\usepackage{pstricks}
\usepackage{graphicx}
\title{A Refinement of the Eulerian numbers, and the Joint Distribution of $\pi(1)$ and $\des(\pi)$ in $S_n$}
\author{Mark Conger\thanks{Department of Mathematics, University of Michigan, 525 East University Avenue, Ann Arbor, MI 48109-1109, USA.  Email: \tt{mconger@umich.edu}}}
\date{July, 2005}
\begin{document}

\pagestyle{plain}
\pagenumbering{arabic}
\bibliographystyle{alpha}
\parindent=0pt
\parskip=12pt
\input{defs.tex}

\maketitle
\begin{abstract}
Given a permutation $\pi$ chosen uniformly from $S_n$, we explore the joint distribution of
$\pi(1)$ and the number of descents in $\pi$.  We obtain a formula for
the number of permutations with $\des(\pi)=d$ and $\pi(1)=k$, and use it to show
that if $\des(\pi)$ is fixed at $d$, then the expected value of $\pi(1)$ is $d+1$.
We go on to derive generating functions for the joint distribution, show that it is
unimodal if viewed correctly, and show that when $d$ is small the distribution of
$\pi(1)$ among the permutations with $d$ descents is approximately geometric.  Applications
to Stein's method and the Neggers-Stanley problem are presented.
\end{abstract}

\section{Introduction} \label{sec:intro}

Consider $S_n$ to be the set of all bijections from $\set{1,2,\ldots,n}$ to itself.
We will often identify a permutation $\pi$ with the sequence $\pi(1),\pi(2),\ldots,\pi(n)$.  So
for instance if $\pi(1)=k$ and $\pi(n)=\ell$, we say that $\pi$ ``begins with'' $k$ and
``ends with'' $\ell$.

A permutation $\pi$ is said to have a descent at $i$ if $\pi(i)>\pi(i+1)$.  That is to say, if we graph
the points $(i,\pi(i))$ and connect them left to right, descents are the positions at which the
connecting segments have negative slope.  Let $\des(\pi)$ be the number of descents in $\pi$, and define
\begin{equation} \label{eulerdef}
	\euler{n}{d} \assign \sizeof{\set{\pi\in S_n:\des(\pi)=d}}.
\end{equation}
These are known as the Eulerian numbers, and have been widely studied; see, for example, \cite[p. 267]{ConcreteMath}
and \cite{Carlitz}.

Bayer and Diaconis \cite{BayerDiaconis} showed that the probability that a 
particular permutation of a deck of cards occurs after any number of riffle 
shuffles is determined by the number of descents the permutation has.  In
\cite{Shuffle1}, Viswanath and the author began working to generalize
that result to decks containing repeated cards.  At one point we had occasion to
consider the number of permutations of $n$ letters which have $d$ descents and begin with $k$.
That is,
\begin{equation} \label{eq:andkdef}
	\euler{n}{d}_k \assign \sizeof{\set{\pi \in S_n : \mbox{$\des(\pi) = d$ and $\pi(1) = k$}}}.
\end{equation}
The current work is an investigation of the numbers defined in \eqnum{andkdef}.  We derive
a formula in terms of binary coefficients:
\begin{prototheorem}{andkformula}
	If $1 \le k \le n$,
	\[ \euler{n}{d}_k = \sum_{j\ge 0} (-1)^{d-j} \binom{n}{d-j} j^{k-1}(j+1)^{n-k} \]
	where $0^0$ is interpreted as 1
\end{prototheorem}
which is similar to a well-known formula for the Eulerian numbers.  We use the formula to understand how
the two statistics $\des(\pi)$ and $\pi(1)$ interact.

If we are constructing a permutation with $d$
descents from left to right, and $d$ is small, a conservative strategy would seem to be to start with a
low number, since starting with a high number means we will use up one of our descents near the beginning
of the permutation.  So in other words, we expect that if $d$ is small then there are more permutations
with $d$ descents starting with low numbers than starting with high numbers.  Similarly, if $d$ is
close to $n$, our intuition is that that starting with a high number leaves us more possibilities later on.
This intuition turns into a surprisingly simple result:
\begin{prototheorem}{edpi1}
If $\pi$ is chosen uniformly from among those permutations of $n$ that have $d$ descents,
the expected value of $\pi(1)$ is $d+1$ and the expected value of $\pi(n)$ is $n-d$.
\end{prototheorem}
And in \thmnum{unimodal} we find, as expected, that the sequence
\[ \euler{n}{d}_1,\euler{n}{d}_2,\ldots,\euler{n}{d}_n \]
is weakly decreasing when $d$ is small and weakly increasing when $d$ is 
large.  Consequently that sequence is an interpolation between its 
endpoints, which are two Eulerian numbers: $\euler{n-1}{d}$ and $\euler{n-1}{d-1}$.
Experimental evidence (see \secnum{remarks}) suggests that it 
is a good interpolation, at least when $d$ is close to $(n-1)/2$, in the 
sense that a normal approximation to the Eulerian numbers also seems to 
provide a good approximation to the refined Eulerian numbers.  However, the
normal approximation is good for neither set when $d$ is small or $d$ is close to
$n$.  \thmnum{unimodal} shows that in those cases the distribution of $\pi(1)$ is
approximately geometric.

The application which led directly to the current work is presented in \secnum{stein}.
Fulman shows in \cite{Fulman1} that certain statistics on permutations, one of which is descents,
are approximately normally distributed.  The main tool he uses is Stein's method,
due to Charles Stein in \cite{Stein}.  The thrust behind the method is to introduce a little
extra randomness to a given random variable to get a new one.  If certain symmetries are present, the
result is an ``exchangeable pair'' of random variables, meaning, essentially, that the Markov
process which takes one to the other is reversible.  Then Stein's theorems (and more recent
refinements of them) can be applied to bound the distance between the original variable's
distribution and the standard normal distribution.

Fulman uses a ``random to end'' operation to add randomness to permutations.  That is,
he starts with a uniformly distributed permutation $\pi$ and sets
\[ \pi' = (I,I+1,\ldots,n)\pi \]
where $I$ is selected uniformly from $\set{1,2,\ldots,n}$.  While $(\pi,\pi')$ is not an exchangeable
pair, it turns out that $(\des(\pi),\des(\pi'))$ is, and this leads to a central limit theorem
for descents, and for a whole class of statistics.

We tried a different method of adding randomness to $\pi$, namely, following $\pi$ by
a uniformly selected transposition.  That calculation (which is presented in \secnum{stein})
led directly to \thmnum{edpi1}.

The Neggers-Stanley Conjecture, now proved false in general 
(\cite{Branden1,Stembridge1}), was that the generating function for descents 
among the linear extensions of any poset has only real zeroes.  Since a 
function with positive coefficients can have no positive zeroes, any combinatorial generating
function with all real zeroes can be written in the form
\[ a(x+c_1)(x+c_2)\cdots(x+c_n) \]
for non-negative constants $a,c_1,c_2,\ldots,c_n$.  The implication, then, is that
if $D$ is the number of descents in a uniformly selected linear extension of a poset for
which the Neggers-Stanley conjecture is true, then $D$ can be written as the sum of
independent Bernoulli variables.

In \secnum{genfuncs} we present several generating functions for the refined Eulerian numbers.
The set of permutations of $n$ which begin with $k$ is the same as the set of linear extensions
of the poset defined on $\set{1,2,\ldots,n}$ by $k<a$ for all $a$ other than $k$.  So we can find the
Neggers-Stanley generating function for this poset explicitly, and we show that it does indeed
have only real zeroes.  We go on to show that several similar posets also satisfy the conjecture.
(All of the posets considered were known to satisfy the conjecture
by theorems of Simion \cite{Simion} and Wagner \cite{Wagner}.)

\section{Basic Properties} \label{sec:basic}

If $\pi(1) = 1$, then $\pi(1)$ is certainly less than $\pi(2)$, so all descents are
among the final $n-1$ numbers.  And if $\pi(1) = n$, there is certain to be a descent
between $\pi(1)$ and $\pi(2)$.  So we know some boundary values:
\begin{equation} \euler{n}{d}_1 = \euler{n-1}{d} \andd \euler{n}{d}_n = \euler{n-1}{d-1} \label{eq:boundries} \end{equation}
for $n > 1$.  Also, it is immediate that
\begin{align}
	\sum_d \euler{n}{d}_k& = (n-1)! \\
	\sum_k \euler{n}{d}_k& = \euler{n}{d}.
\end{align}
Let $\rho \in S_n$ be the reversal permutation: $\rho(i) = n+1-i$.  Then $\rho\pi$ is the same as
$\pi$ but with $i$ replaced by $n+1-i$ everywhere.  As a result, $\rho\pi$ has a descent wherever
$\pi$ has an ascent, and an ascent wherever $\pi$ has a descent.  So $\des(\rho\pi) = n-1-\des(\pi)$.  Since $\pi \mapsto \rho\pi$
is a bijection from $S_n$ to itself, it follows that
\begin{equation} \euler{n}{d} = \euler{n}{n-1-d}. \label{eq:eulersym} \end{equation}
Note we could have obtained the same result from the map $\pi \mapsto \pi\rho$, since reversing
$\pi$ changes ascents to descents and also reflects their positions about the center.

Let
\begin{equation} \euler{n}{d}^k \assign \sizeof{\set{\pi \in S_n : \mbox{$\des(\pi) = d$ and $\pi(n) = k$}}}. \label{eq:supdef} \end{equation}
Both transformations yield symmetric identities for the refined Eulerian numbers.  If
\[ \pi(1)=k \andd \des(\pi)=d \]
then
\begin{align*}
	\rho\pi(1)=n+1-k &\andd \des(\rho\pi)=n-1-d \\
	\pi\rho(n)=k &\andd \des(\pi\rho)=n-1-d \\
	\rho\pi\rho(n)=n+1-k &\andd \des(\rho\pi\rho)=d
\end{align*}
from which it follows that
\begin{equation} \label{eq:symmetry}
	\euler{n}{d}_k = \euler{n}{n-1-d}_{n+1-k} = \euler{n}{n-1-d}^k = \euler{n}{d}^{n+1-k}.
\end{equation}

\section{Recurrences} \label{sec:recurrences}

Assume $n>1$.  Let
\begin{align*}
	T_k &\assign \set{\pi \in S_n: \mbox{$\pi(1)=k$ and $\des(\pi)=d$}} \\
	T_{k,\ell} &\assign \set{\pi \in S_n: \mbox{$\pi(1)=k$, $\pi(2)=\ell$, and $\des(\pi)=d$}}
\end{align*}
and let $\pi \in T_{k,\ell}$.
If $\ell<k$, then there is a descent between $\pi(1)$ and $\pi(2)$, so there must be $d-1$ descents
in the ``tail'', $\pi(2),\pi(3),\ldots,\pi(n)$.  The tail begins with $\ell$, which is the $\ell$th largest value in the tail,
so we must have
\[ \sizeof{T_{k,\ell}} = \euler{n-1}{d-1}_\ell \]
when $\ell<k$.  Likewise, if $\ell>k$, there is no descent between $\pi(1)$ and $\pi(2)$, so there must be
$d$ descents in the tail.  This time $\ell$ is the $(\ell-1)$st largest value in the tail, so
\[ \sizeof{T_{k,\ell}} = \euler{n-1}{d}_{\ell-1} \]
when $\ell>k$.  Of course $T_k$ is the disjoint union of the $T_{k,l}$, so
\[ \euler{n}{d}_k = \sizeof{T_k} = \sum_\ell \sizeof{T_{k,\ell}} = \sum_{\ell<k} \euler{n-1}{d-1}_\ell + \sum_{\ell>k} \euler{n-1}{d}_{\ell-1} \]
or, more succinctly,
\begin{equation} \label{eq:rec1} 
	\euler{n}{d}_k =  \sum_{\ell=1}^{n-1} \euler{n-1}{d-[\ell<k]}_\ell
\end{equation}
where the bracket notation follows \cite{ConcreteMath}: $[A]$ is 1 if $A$ is true and 0 if $A$ is false.
(Knuth refers to this as Iverson notation in \cite{KnuthNotation}, and traces its origin to \cite{APL}.)
Note \eqnum{rec1} fails when $k<1$ or $k>n$, in which case ${\textstyle \euler{n}{d}_k} = 0$.


Now suppose $1 \le k \le n-1$ and $\pi\in S_n$ begins with $k$.  Swapping $k$ with $k+1$
in the sequence $\pi(1),\pi(2),\ldots,\pi(n)$ preserves descents for most $\pi$; the only
exception is when $\pi(2)=k+1$, in which case a new descent is created.  If we eliminate
that case, the swap map is a bijection from $T_k \setminus T_{k,k+1}$ to $T_{k+1} \setminus T_{k+1,k}$,
as those sets are defined above.  Substituting sizes for sets, we have
\begin{equation} \label{eq:rec2}
	\euler{n}{d}_k - \euler{n-1}{d}_k = \euler{n}{d}_{k+1} - \euler{n-1}{d-1}_k.
\end{equation}
\eqnum{rec2} is valid as long as $k \ne 0$ and $k \ne n$.  (If $k < 0$ or $k > n$, all terms are 0.)

A well-know recurrence for $\euler{n}{d}$ comes from considering
what happens when you insert $n$ into an element of $S_{n-1}$:
\begin{equation} \label{eq:eulerrec}
	\euler{n}{d} = (n-d)\euler{n-1}{d-1} + (d+1)\euler{n-1}{d}
\end{equation}
We can get a similar recurrence for the refined Eulerian numbers by 
considering what happens when you insert $n$ into an element of $S_{n-1}$ which begins with $k$:
\begin{equation} \label{eq:rec3}
	\euler{n}{d}_k = (n-d-1)\euler{n-1}{d-1}_k + (d+1)\euler{n-1}{d}_k.
\end{equation}
In other words, one way to get an element of $S_n$ which begins with $k$ and has $d$ descents
is to take an element of $S_{n-1}$ which begins with $k$ and has $d$ descents,
and insert $n$ at a descent or at the end ($d+1$ choices).  The other way is to start with an element
of $S_{n-1}$ which begins with $k$ and has $d-1$ descents, and insert $n$ at an ascent
($n-d-1$ choices).
\eqnum{rec3} fails when $k=n$, since a permutation of $S_{n-1}$ cannot begin with $n$.
It is valid for all other values of $k$.

\section{Formulas and Moments} \label{sec:formulas}

There is an explicit formula for the Eulerian numbers in terms of binomial coefficients:
\begin{equation} \label{eq:eulerformula}
	\euler{n}{d} = \sum_{j\ge 0} (-1)^{d-j} \binom{n+1}{d-j}(j+1)^n.
\end{equation}
See for example, \cite[p. 269]{ConcreteMath}.  [Aside: \eqnum{eulerformula} follows from 
\eqnum{eulerrec}, which means that it is valid for all values of $d$, even 
if $d<0$ or $d\ge n$].  So we have
\begin{align}
	\euler{n}{d}_1 &= \euler{n-1}{d} = \sum_{j\ge 0} (-1)^{d-j} \binom{n}{d-j}(j+1)^{n-1} \label{eq:formulaone} \\
	\euler{n}{d}_n &= \euler{n-1}{d-1} = \sum_{j\ge 0} (-1)^{d-1-j} \binom{n}{d-1-j}(j+1)^{n-1} \label{eq:formulan}
	               = \sum_{j\ge 0} (-1)^{d-j} \binom{n}{d-j} j^{n-1}.
\end{align}
These suggest a formula for $\euler{n}{d}_k$:

\begin{theorem} \label{thm:andkformula}
	If $1 \le k \le n$,
	\begin{equation} \label{eq:andkformula}
		\euler{n}{d}_k = \sum_{j\ge 0} (-1)^{d-j} \binom{n}{d-j} j^{k-1}(j+1)^{n-k}
	\end{equation}
	where $0^0$ is interpreted as 1.
\end{theorem}

\proof Fix $k\ge 1$.  The theorem is true for $n=k$ by \eqnum{formulan}.  Suppose it is true for some $n\ge k$.
Then
\begin{align*}
	\euler{n+1}{d}_k &= (n-d)\euler{n}{d-1}_k + (d+1)\euler{n}{d}_k \\
	                 &= \sum_{j\ge 0} (-1)^{d-j}\bracket{-(n-d)\binom{n}{d-1-j} + (d+1)\binom{n}{d-j}} j^{k-1}(j+1)^{n-k}.
\end{align*}
The quantity in brackets reduces to $(j+1)\binom{n+1}{d-j}$, so the theorem is true for all $n\ge k$ by induction.
\eop

Note that we assumed nothing about $d$; \eqnum{andkformula} is valid even if $d<0$ or $d\ge n$.

From \eqnum{andkformula} we can deduce a formula for the $m$th ``rising moment'' of $\pi(1)$ when
$\des(\pi)$ is fixed.  Assume $\pi$ is chosen uniformly from $S_n$, and let
\begin{equation} \mu_m \assign \ex^{\des(\pi)=d} \pi(1)^{\overline{m}} \label{eq:mudef} \end{equation}
where $x^{\overline{m}} = x(x+1)(x+2)\cdots(x+m-1)$.

\begin{lemma}
\begin{equation}
	\euler{n}{d}\mu_m = m! \sum_{j \ge 0} (-1)^{d-j} \binom{n}{d-j} \sum_{\ell=0}^{n-1} \binom{m+n}{\ell} j^{\ell}.
	\label{eq:momentformula}
\end{equation}
\end{lemma}
\proof From \eqnum{andkformula},
\begin{align*}
	\euler{n}{d}\mu_m &= \sum_{k=1}^n k^{\overline{m}} \euler{n}{d}_k
	                   = \sum_{k=1}^n \frac{(k+m-1)!}{(k-1)!} \sum_{j\ge 0} (-1)^{d-j} \binom{n}{d-j} j^{k-1}(j+1)^{n-k} \\
	                  &= m! \sum_{j\ge 0} (-1)^{d-j} \binom{n}{d-j} \sum_{r=0}^{n-1} \binom{r+m}{r} j^r(j+1)^{n-1-r}
\end{align*}
(the last by setting $r=k-1$).
But $(j+1)^{n-1-r} = \sum_{s=0}^{n-1-r} \binom{n-1-r}{s} j^s$.  So let $\ell=r+s$ and we have
\begin{equation} \label{eq:muformcomplicated}
	\euler{n}{d}\mu_m = m! \sum_{j\ge 0} (-1)^{d-j} \binom{n}{d-j} \sum_{\ell=0}^{n-1} j^\ell \sum_{r=0}^\ell \binom{r+m}{r} \binom{n-1-r}{\ell-r}.
\end{equation}
\newcommand{\xtick}[1]{\psline[linestyle=solid,linewidth=.5pt](#1,-.25)(#1,.25)}
\newcommand{\ytick}[1]{\psline[linestyle=solid,linewidth=.5pt](-.25,#1)(.25,#1)}
\newcommand{\xlabel}[3]{\xtick{#1}\rput[#3](#1,-.75){#2}}
\newcommand{\xlabels}[2]{\xlabel{#1}{#2}{B}}
\newcommand{\xlabelse}[2]{\xlabel{#1}{#2}{Bl}}
\newcommand{\ylabelw}[2]{\ytick{#1}\uput[l](-.25,#1){#2}}
\begin{figure}
	\centering
	\psset{unit=20pt}
	\begin{pspicture}(-1,-1)(11,9)
		\psline[linewidth=1pt]{->}(0,0)(10,0)
		\uput[r](10,0){$x$}
		\multips(0,1)(0,1){7}{\psline[linestyle=dotted,linewidth=.5pt](0,0)(9,0)}
		\xlabels{0}{0}
		\xlabels{5}{$m$}
		\xlabelse{6}{$m+1$}
		\xlabelse{9}{$m+n-\ell$}

		\psline[linewidth=1pt]{->}(0,0)(0,8)
		\uput[u](0,8){$y$}
		\multips(1,0)(1,0){9}{\psline[linestyle=dotted,linewidth=.5pt](0,0)(0,7)}
		\ylabelw{0}{0}
		\ylabelw{4}{$r$}
		\ylabelw{7}{$\ell$}

		\psline[linewidth=2pt](0,0)(1,0)(1,2)(2,2)(2,3)(5,3)(5,4)(7,4)(7,6)(9,6)(9,7)
		\pscircle*(0,0){.15}
		\pscircle*(5,4){.15}
		\pscircle*(6,4){.15}
		\pscircle*(9,7){.15}
		\psline[linewidth=1pt,linestyle=dashed](5.5,-1)(5.5,8)
		\uput[r](5.5,8){$x=m+\frac12$}
	\end{pspicture}
	\caption{A north-east lattice path from $(0,0)$ to $(m+n-\ell,\ell)$.  (All edges are either north or east.)}
	\label{fig:latticepath}
\end{figure}
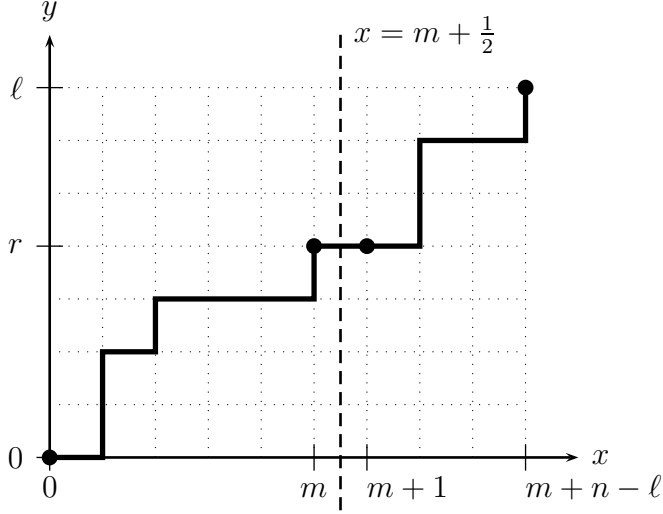
Let $\phi$ be a north/east lattice path from $(0,0)$ to $(m+n-\ell,\ell)$ (see Figure~\ref{fig:latticepath}).
The number of such paths is $\binom{m+n}{\ell}$.  If $r$ is the height at which $\phi$ crosses the line $x=m+\frac12$,
then $\phi$ consists of a path from $(0,0)$ to $(m,r)$, a horizontal segment, and a path from $(m+1,r)$
to $(m+n-\ell,\ell)$.  Counting the possibilities for the parts yields the identity
\begin{equation} \label{eq:lpidentity}
	\sum_{r=0}^\ell \binom{r+m}{r} \binom{n-1-r}{\ell-r} = \binom{m+n}{\ell}.
\end{equation}
Substituting \eqnum{lpidentity} into \eqnum{muformcomplicated} yields the desired result.
\eop

Note that the last sum in \eqnum{momentformula} is a truncated binomial expansion of $(j+1)^{m+n}$.

\begin{theorem} \label{thm:edpi1}
If $\pi$ is chosen uniformly from among those permutations of $n$ that have $d$ descents,
the expected value of $\pi(1)$ is $d+1$ and the expected value of $\pi(n)$ is $n-d$.
\end{theorem}

\proof The expected value of $\pi(1)$ is $\mu_1$, and
\begin{align*}
	\euler{n}{d} \mu_1& = \sum_{j \ge 0} (-1)^{d-j} \binom{n}{d-j} \sum_{\ell=0}^{n-1} \binom{n+1}{\ell} j^{\ell} \\
	& = \sum_{j \ge 0} (-1)^{d-j} \binom{n}{d-j} \paren{(j+1)^{n+1} - j^{n+1} - (n+1)j^n} \\
	& = \sum_{j \ge 0} (-1)^{d-j} \binom{n}{d-j} (j+1)^{n+1} - \sum_{i \ge 0} (-1)^{d-i} \binom{n}{d-i}(n+1+i)i^n.
	\intertext{The term for $i=0$ is 0, so let $j=i-1$ and combine}
	& = \sum_{j \ge 0} (-1)^{d-j} (j+1)^n \bracket{\binom{n}{d-j}(j+1) + \binom{n}{d-j-1}(n+j+2)}.
	\intertext{The quantity in brackets simplifies to $(d+1)\binom{n+1}{d-j}$, so}
	\euler{n}{d} \mu_1& = (d+1) \sum_{j \ge 0} (-1)^{d-j} (j+1)^n \binom{n+1}{d-j} = (d+1) \euler{n}{d}.
\end{align*}
Therefore
\[ \mu_1 = \ex^{\des(\pi)=d} \pi(1) = d+1. \]
For the second part,
\begin{align*}
	\ex^{\des(\pi)=d} \pi(n) &= \frac1{n!} \sum_k k\euler{n}{d}^k = \frac1{n!} \sum_k k\euler{n}{n-1-d}_k \\
	                         &= \ex^{\des(\pi)=n-1-d} \pi(1) = n-d.
\end{align*}
\eop

\section{Application Using Stein's Method} \label{sec:stein}

Charles Stein developed a method for showing that the distribution of a 
random variable $W$ which meets certain criteria is approximately standard normal.  His 
technique has come to be known as Stein's method; see \cite{Stein} or \cite{SteinApplications}
for more explanation than can be given here.

In its most straightforward form, Stein's method requires finding a ``companion''
random variable $W^*$ such that $(W,W^*)$ is an exchangeable pair, meaning that
\begin{equation} \label{eq:exchangeable}
	\pr(W=w,W^*=w^*) = \pr(W=w^*,W^*=w)
\end{equation}
for all values of $w$ and $w^*$.  If we can find such a $W^*$ and if, in addition,
there is a $\lambda$ between 0 and 1 such that
\begin{equation} \label{eq:lambda}
	\ex^W W^* = (1-\lambda)W
\end{equation}
(that is, the expected value of $W^*$ when $W$ is fixed at some value is $1-\lambda$
times that value), then we may apply Stein's method.

We are interested in showing that if $\pi$ is chosen uniformly from $S_n$, then
the random variable
$D = \des(\pi)$ is approximately normal.  This has been proven before, and in more generality; see
\cite{Fulman1} for references.  We will demonstrate the set-up for Stein's method---that is, finding
a companion variable and showing that it satisfies \eqnum{exchangeable} and \eqnum{lambda}.  From there,
applying the method would proceed as in \cite{Fulman1}.

Often the companion variable in Stein's method is defined by adding a little
bit of randomness to the variable we are interested in.  In this case, let
$\tau$ be selected uniformly from among the transpositions in $S_n$, independently
of $\pi$.  Then $\tau\pi$ is uniformly distributed over $S_n$, and for any $u,v\in S_n$,
\begin{align*}
	\pr(\pi=u,\tau\pi=v) &= \pr(\pi=u,\tau=vu^{-1}) = \pr(\pi=u)\pr(\tau=vu^{-1}) \\
	\pr(\pi=v,\tau\pi=u) &= \pr(\pi=v,\tau=uv^{-1}) = \pr(\pi=v)\pr(\tau=uv^{-1}).
\end{align*}
Both right-hand sides are $(n!)^{-1}\binom{n}{2}^{-1}$ if $vu^{-1}$ is a transposition and 0 otherwise,
so $(\pi,\tau\pi)$ is an exchangeable pair.

Let $D^* \assign \des(\tau\pi)$.  Since $(\pi,\tau\pi)$ is an exchangeable pair, $(F(\pi),F(\tau\pi))$
is exchangeable for any function $F$.  So $(D,D^*)$ is exchangeable.  For $1 \le i \le n-1$ let
\[ D_i = [\pi(i) > \pi(i+1)] \andd D_i^* = [\tau\pi(i) > \tau\pi(i+1)] \]
be Bernoulli random variables; then $D = \sum_{i=1}^{n-1} D_i$ and $D^* = \sum_{i=1}^{n-1} D_i^*$.

Fix $\pi$ and $i$ and let $a=\pi(i)$, $b=\pi(i+1)$.  If $a<b$, the only ways for $\tau\pi(i)$
to be bigger than $\tau\pi(i+1)$ are if $\tau$ swaps $a$ with something bigger than $b$ ($n-b$ ways),
if $\tau$ swaps $b$ with something smaller than $a$ ($a-1$ ways), or if $\tau$ swaps $a$ with $b$.
So
\[ \ex^{D_i=0} (D_i^*-D_i) = \pr(D_i^*=1|D_i=0) = \frac{n+\pi(i)-\pi(i+1)}{\binom{n}{2}} \]
and similarly if $a>b$,
\[ \ex^{D_i=1} (D_i^*-D_i) = -\pr(D_i^*=0|D_i=1) = -\frac{n+\pi(i+1)-\pi(i)}{\binom{n}{2}}. \]
So in general
\[ \ex^{D_i} (D_i^*-D_i) = \frac{\pi(i)-\pi(i+1)}{\binom{n}{2}} + \frac{2(1-2D_i)}{n-1}. \]
Summing now over $i$ causes the $\pi(i)$ terms to telescope:
\[ \ex^{\pi} (D^*-D) = \sum_{i=1}^{n-1} \ex^{\pi} (D_i^*-D_i) = \frac{\pi(1)-\pi(n)}{\binom{n}{2}} + 2 - \frac{4D}{n-1} \]
which allows us to apply \thmnum{edpi1}:
\begin{align*}
	\ex^D (D^*-D) &= \ex^D \ex^{\pi} (D^*-D)
	               = \frac{\ex^D \pi(1) - \ex^D \pi(n)}{\binom{n}{2}} + 2 - \frac{4D}{n-1} \\
	              &= \frac{2}{n(n-1)}((D+1)-(n-D)) + 2 - \frac{4D}{n-1}
	               = \frac{2(n-1)-4D}{n}.
\end{align*}
The mean and variance of $\des(\pi)$ are
$\mu \assign (n-1)/2$ and $\sigma^2 \assign (n+1)/12$ respectively, so the variables
\[ W \assign \frac{\des(\pi)-\mu}{\sigma} \andd W^* \assign \frac{\des(\tau\pi)-\mu}{\sigma} \]
have mean 0 and variance 1.  Then $(W,W^*)$ is an exchangeable pair and
\[ \ex^{W=w}(W^*-W) = \ex^{D=\sigma w + \mu} \paren{\frac{D^*-\mu}{\sigma}-\frac{D-\mu}{\sigma}} = \frac{1}{\sigma} \ex^{D=\sigma w + \mu} (D^*-D) \]
which is to say
\[ \ex^W (W^*-W) = \frac{2(n-1)-4(\sigma W + \mu)}{\sigma n} = -\frac4n W. \]
So if $W^*$ is obtained using the ``random transposition'' method described here,
$(W,W^*)$ will be an exchangeable pair satisfying \eqnum{lambda} with $\lambda=4/n$.  One can now
proceed with Stein's method and show that $W$ is close to being a standard normal random variable.

\section{Generating Functions} \label{sec:genfuncs}

It follows from \eqnum{eulerformula} that
\[ a_n(x) \assign \sum_d \euler{n}{d} x^{d+1} = (1-x)^{n+1} \sum_{j \ge 0} j^n x^j \]
and therefore that
\begin{align*}
	A(x,z) :&= \sum_{n \ge 1} a_n(x) z^n/n!
	         = \sum_{n \ge 1} (1-x)^{n+1} \sum_{j \ge 0} j^n x^j z^n / n! \\
	        &= (1-x) \sum_{j \ge 0} x^j \sum_{n \ge 1} \paren{j(1-x)z}^n / n!
	         = (1-x) \sum_{j \ge 0} x^j \paren{e^{j(1-x)z} - 1} \\
	        &= (1-x) \paren{\frac{1}{1-xe^{(1-x)z}} - \frac{1}{1-x}}
	         = \frac{1-x}{1-xe^{(1-x)z}} - 1 \\
	        &= \frac{xe^{-(1-x)z} - x}{x - e^{-(1-x)z}}.
\end{align*}
There is some disagreement in the literature about what $a_0$ should be.  We have avoided the
problem by not including it in the sum.

There are various ways to define generating functions for the $\euler{n}{d}_k$, depending on which
variables are kept constant.

\begin{theorem} \label{thm:gfall}
\begin{equation} \label{eq:gfall}
	\sum_{n,d,k} \euler{n}{d}_k x^d y^k z^n/n! = \frac{1}{\theta} \int_\theta^{\theta^y} \frac{dt}{x - t^{1-1/y}}
\end{equation}
where $\theta = \exp\left\{\paren{\frac{1-x}{y^{-1}-1}}z\right\}$.
\end{theorem}

\proof Let $B(x,y,z)$ be the left-hand side of \eqnum{gfall}.  Note the sum is over all integers $n$, $d$, and $k$.
So
\begin{eqnarray} \label{eq:gfdea}
	\lefteqn{(y^{-1}-1)\del{B}{z} + (1-x)B = } \\
	& & \sum_{n,d,k}
		\bracket{\euler{n+1}{d}_{k+1} - \euler{n+1}{d}_k + \euler{n}{d}_k - \euler{n}{d-1}_k}
	x^d y^k z^n/n! \nonumber
\end{eqnarray}
Let $S(n,d,k)$ be the bracketed quantity.  It is clearly 0 if $n<0$, and if $n=0$,
\[ S(0,d,k) = \begin{piecewise}
	\pwif{1}{d=0,k=0} \\
	\pwif{-1}{d=0,k=1} \\
	\pwot{0}
\end{piecewise} \]
so $n=0$ contributes $1-y$ to the sum on the right-hand side of \eqnum{gfdea}.  If $n\ge 1$, then
by \eqnum{rec2}, $S(n,d,k)$ is 0 unless $k=0$ or $k=n+1$, in which case
\[ S(n,d,0) = \euler{n+1}{d}_1 = \euler{n}{d} \andd S(n,d,n+1) = -\euler{n+1}{d}_{n+1} = -\euler{n}{d-1}. \]
Therefore
\begin{align*}
	(y^{-1}-1)\del{B}{z} + (1-x)B
	&= 1-y + \sum_{n \ge 1,d} \euler{n}{d} x^d z^n / n!
	- \sum_{n \ge 1,d} \euler{n}{d-1} x^d y^{n+1} z^n / n! \\
	&= 1-y + x^{-1}A(x,z) - yA(x,yz) \\
	&= \bracket{1 + x^{-1}\frac{xe^{-(1-x)z}-x}{x-e^{-(1-x)z}}}
	 - y \bracket{1 + \frac{xe^{-(1-x)yz}-x}{x-e^{-(1-x)yz}}} \\
	&= (1-x) \bracket{ \frac{ye^{-(1-x)yz}}{x-e^{-(1-x)yz}} - \frac{1}{x-e^{-(1-x)z}} }.
\end{align*}
Let $\alpha = \frac{1-x}{y^{-1}-1}$. Then $\theta$, as defined in the theorem, is $e^{\alpha z}$.
Dividing by $y^{-1}-1$ and multiplying through by $\theta$ gives
\[ \theta \del{B}{z} + \alpha\theta B = \alpha\theta \bracket{ \frac{ye^{-(1-x)yz}}{x-e^{-(1-x)yz}} - \frac{1}{x-e^{-(1-x)z}} } \]
which is to say that
\[ \del{}{z}\paren{\theta B} = \alpha \bracket{ \frac{y\theta^y}{x-\theta^{y-1}} - \frac{\theta}{x-\theta^{1-1/y}} }. \]
Differentiating the integral on the right-hand side of \eqnum{gfall},
\begin{align*} \del{}{z} \int_\theta^{\theta^y} \frac{dt}{x - t^{1-1/y}}
	&= \del{\theta^y}{z} \bracket{\frac{1}{x-\paren{\theta^y}^{1-1/y}}} - \del{\theta}{z} \bracket{\frac{1}{x-\theta^{1-1/y}}} \\
	&= \frac{\alpha y \theta^y}{x-\theta^{y-1}} - \frac{\alpha\theta}{x-\theta^{1-1/y}} = \del{}{z}\paren{\theta B}.
\end{align*}
Since $\theta B$ and the integral have the same derivative with respect to $z$, and they both vanish when
$z=0$, they are equal. \eop

Here are three more generating functions.  They can all be found
by plugging in \eqnum{andkformula} and switching summation signs.
\begin{align}
	\sum_d &\euler{n}{d}_k x^d = (1-x)^n \sum_{j\ge 0} j^{k-1}(j+1)^{n-k}x^j \label{eq:gfnk} \\
	\sum_k &\euler{n}{d}_k y^k = y \sum_{j\ge 0} (-1)^{d-j} \binom{n}{d-j} \frac{(j+1)^n - (jy)^n}{j+1-jy} \label{eq:gfnd} \\
	\sum_{d,k} &\euler{n}{d}_k x^d y^k = (1-x)^n y \sum_{j \ge 0} \frac{(j+1)^n - (jy)^n}{j+1-jy} x^j \label{eq:gfn}
\end{align}

We can now prove a special case of the Neggers-Stanley conjecture.
Define the descent polynomial of $A \subset S_n$ to be
\[ F_A(x) = \sum_{\pi\in A} x^{\des(\pi)}. \]
Let $P$ be a poset of $n$ elements with labels $1,2,\ldots,n$.  A linear extension of $P$
is an ordering of $1,2,\ldots,n$ which preserves the ordering of $P$; that is, a $\pi\in S_n$
which is such that if $i <_P j$ then $i$ appears before $j$ in the list $\pi(1),\pi(2),\ldots,\pi(n)$.
If $\mathcal{L}(P)$ denotes the set of linear extentions of $P$, then Neggers and Stanley \cite[p. 311]{StanleyConj}
conjectured that for any poset, every zero of $F_{\mathcal{L}(P)}$ is real.

The conjecture has been shown to be false in general \cite{Branden1,Stembridge1}.  But we can prove it is
true in a certain special case.

\begin{theorem} \label{thm:neggers1}
	If $P_{n,k}$ is the poset with Hasse diagram
	\begin{center}
	\psset{xunit=.6in,yunit=.5in}
	\begin{pspicture}(-4,-.1)(4,1.1)
		\psline[linewidth=.5pt](0,0)(-3.5,1)
		\psline[linewidth=.5pt](0,0)(-2.5,1)
		\psline[linewidth=.5pt](0,0)(-0.5,1)
		\psline[linewidth=.5pt](0,0)(0.5,1)
		\psline[linewidth=.5pt](0,0)(1.5,1)
		\psline[linewidth=.5pt](0,0)(3.5,1)
		\uput[d](0,0){$k$}
		\uput[u](-3.5,1){$1$}
		\uput[u](-2.5,1){$2$}
		\uput[u](-1.5,1){$\cdots$}
		\uput[u](-0.5,1){$k-1$}
		\uput[u](0.5,1){$k+1$}
		\uput[u](1.5,1){$k+2$}
		\uput[u](2.5,1){$\cdots$}
		\uput[u](3.5,1){$n$}
		\psdots[dotstyle=*](0,0)(-3.5,1)(-2.5,1)(-0.5,1)(0.5,1)(1.5,1)(3.5,1)
	\end{pspicture}
	\end{center}
	then $F_{\mathcal{L}(P_{n,k})}$ has only distinct real roots.
\end{theorem}

\proof For $u,v \ge 0$ let
\[
	c_{u,v} \assign \sum_d \euler{u+v+1}{d}_{u+1}x^d = \sum_{\substack{
		\pi\in S_{u+v+1} \\
		\pi(1)=u+1
	}} x^{\des(\pi)}.
\]
Then setting $u=k-1$, $v=n-k$ yields the polynomial in question.  If $v=0$,
$c_{u,v}$ counts the reversal permutation $\rho$, which has $(u+v+1)-1 = u$
descents.  Otherwise, if $v>0$, $c_{u,v}$ doesn't count $\rho$ but it does count
the permutation
\[ u+1,u,u-1,\ldots,1,u+v+1,u+v,\ldots,u+2 \]
which has u+v-1 descents.  So
\[ \deg(c_{u,v}) = \begin{piecewise}
	\pwif{u}{v=0} \\
	\pwif{u+v-1}{v>0.}
\end{piecewise} \]
Similarly, if $u=0$, $c_{u,v}$ counts the identity permutation, which has no descents.
Otherwise it doesn't count the identity but it does count
\[ u+1,u+2,\ldots,u+v+1,1,2,\ldots,u \]
which has 1 descent.  So $x \nmid c_{0,v}(x)$ and if $u>0$,
$x \mid c_{u,v}(x)$ but $x^2 \nmid c_{u,v}(x)$.
Now let
\[ h_{u,v} \assign \frac{c_{u,v}}{(1-x)^{u+v+1}}. \]
Note that $c_{u,v}(1) = \sizeof{\set{\pi\in S_{u+v+1}:\pi(1)=u+1}} = (u+v)!$, so
$c_{u,v}$ does not have a zero at $x=1$.  Therefore $h_{u,v}$ has exactly the same zeroes
as $c_{u,v}$, plus a pole at $x=1$.
By \eqnum{gfnk},
\[ h_{u,v}(x) = \sum_{j\ge 0} j^u (j+1)^v x^j. \]
If $D$ represents differentiation with respect to $x$, we have
\[ (xD) h_{u,v}(x) = h_{u+1,v}(x) \andd (Dx) h_{u,v}(x) = h_{u,v+1}(x) \]
and so
\[ h_{0,v}(x) = (Dx)^v h_{0,0}(x) \andd h_{u,v}(x) = (xD)^u h_{0,v}(x). \]
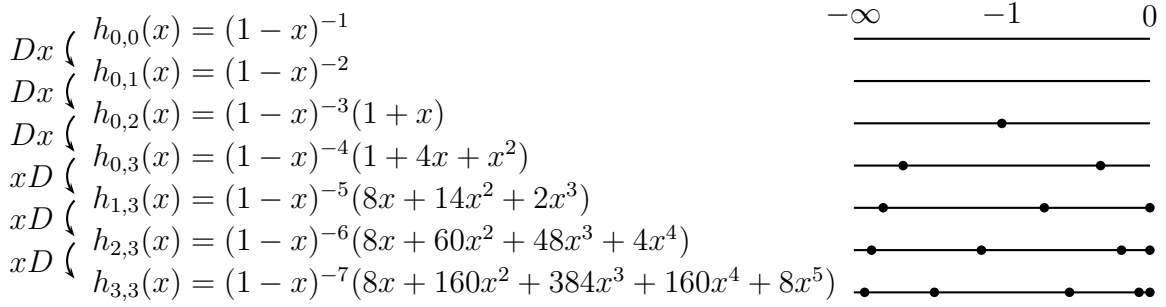
\begin{figure}
	\centering
	\psset{unit=16pt}
	\begin{pspicture}(3,-7)(30,0)
		\rput[Bl]{0}(5, 0){$h_{0,0}(x)=(1-x)^{-1}$}
		\rput[Bl]{0}(5,-1){$h_{0,1}(x)=(1-x)^{-2}$}
		\rput[Bl]{0}(5,-2){$h_{0,2}(x)=(1-x)^{-3}(1+x)$}
		\rput[Bl]{0}(5,-3){$h_{0,3}(x)=(1-x)^{-4}(1+4x+x^2)$}
		\rput[Bl]{0}(5,-4){$h_{1,3}(x)=(1-x)^{-5}(8x+14x^2+2x^3)$}
		\rput[Bl]{0}(5,-5){$h_{2,3}(x)=(1-x)^{-6}(8x+60x^2+48x^3+4x^4)$}
		\rput[Bl]{0}(5,-6){$h_{3,3}(x)=(1-x)^{-7}(8x+160x^2+384x^3+160x^4+8x^5)$}
		\psarc{->}(5,-0.25){.6}{135}{225}
		\psarc{->}(5,-1.25){.6}{135}{225}
		\psarc{->}(5,-2.25){.6}{135}{225}
		\psarc{->}(5,-3.25){.6}{135}{225}
		\psarc{->}(5,-4.25){.6}{135}{225}
		\psarc{->}(5,-5.25){.6}{135}{225}
		\uput[180](4.4,-0.25){$Dx$}
		\uput[180](4.4,-1.25){$Dx$}
		\uput[180](4.4,-2.25){$Dx$}
		\uput[180](4.4,-3.25){$xD$}
		\uput[180](4.4,-4.25){$xD$}
		\uput[180](4.4,-5.25){$xD$}
		\uput[90](23,0){$-\infty$}
		\uput[90](26.5,0){$-1$}
		\uput[90](30,0){$0$}
		\psline(23, 0)(30, 0)
		\psline(23,-1)(30,-1)
		\psline(23,-2)(30,-2)
		\psline(23,-3)(30,-3)
		\psline(23,-4)(30,-4)
		\psline(23,-5)(30,-5)
		\psline(23,-6)(30,-6)
		\psdots[dotstyle=*](26.50000,-2)
		\psdots[dotstyle=*](24.16667,-3)(28.83333,-3)
		\psdots[dotstyle=*](23.69367,-4)(27.50199,-4)(30.00000,-4)
		\psdots[dotstyle=*](23.41905,-5)(26.01482,-5)(29.33017,-5)(30.00000,-5)
		\psdots[dotstyle=*](23.25745,-6)(24.90410,-6)(28.09590,-6)(29.74256,-6)(30.00000,-6)
	\end{pspicture}
	\caption{The construction of $h_{3,3}(x)$ as described in the proof of \thmnum{neggers1}.  The zeroes of each
	function are plotted on the right, using an inverse tangent scale.  Since each function is generated
	from the previous one by applying either the $Dx$ or the $xD$ operator, Rolle's Theorem guarantees
	that the zeroes must interleave.  By a counting argument, all the zeroes of each function must be real.}
	\label{fig:DxxD}
\end{figure}
$h_{0,0}(x) = (1-x)^{-1}$ and $h_{0,1}(x) = (1-x)^{-2}$ both have no zeroes.  Suppose
$v \ge 1$ and $h_{0,v}$ has only distinct real zeroes.  Since $\deg(c_{0,v}) = v-1$ and $x \nmid c_{0,v}(x)$,
$xc_{0,v}(x)$ and $xh_{0,v}(x)$ have $v$ distinct real zeroes.  By Rolle's Theorem,
$(Dx)h_{0,v}$ must have $v-1$ distinct zeroes interlaced between those of $xh_{0,v}(x)$.
Furthermore, the denominator of $xh_{0,v}(x)$ has degree $v+1$, so $xh_{0,v}(x)$ approaches 0 as
$x\rightarrow\infty$.  Therefore its graph must turn back toward the $x$-axis somewhere
to the left of its leftmost zero, at which place there must be another zero of $(Dx)h_{0,v}$.
So we have found $v$ real zeroes of $h_{0,v+1}$, and that accounts for all its zeroes.

Applying the $xD$ operator goes similarly.  Given that $h_{u,v}$ has $d$ distinct
real zeroes, by Rolle's Theorem $Dh_{u,v}(x)$ has $d-1$ interlaced zeroes.  Since the numerator
of $h_{u,v}$ has degree smaller than the denominator, $h_{u,v}$ must turn back toward the axis to the
left of its leftmost zero, which accounts for one more zero of $Dh_{u,v}$.  Finally, $(xD)h_{u,v}$ has
one more zero at 0 (which is distinct from the others since $x^2 \nmid h_{u,v}$ and therefore $x \nmid Dh_{u,v}$).
So we have found $d+1$ real zeroes of $h_{u+1,v}$, and that accounts for all of the zeroes. \eop

\begin{cor} \label{cor:upsidedown} The same can be said for the poset
	\begin{center}
	\psset{xunit=.6in,yunit=.5in}
	\begin{pspicture}(-4,-.1)(4,1.1)
		\psline[linewidth=.5pt](0,1)(-3.5,0)
		\psline[linewidth=.5pt](0,1)(-2.5,0)
		\psline[linewidth=.5pt](0,1)(-0.5,0)
		\psline[linewidth=.5pt](0,1)(0.5,0)
		\psline[linewidth=.5pt](0,1)(1.5,0)
		\psline[linewidth=.5pt](0,1)(3.5,0)
		\uput[u](0,1){$k$}
		\uput[d](-3.5,0){$1$}
		\uput[d](-2.5,0){$2$}
		\uput[d](-1.5,0){$\cdots$}
		\uput[d](-0.5,0){$k-1$}
		\uput[d](0.5,0){$k+1$}
		\uput[d](1.5,0){$k+2$}
		\uput[d](2.5,0){$\cdots$}
		\uput[d](3.5,0){$n$}
		\psdots[dotstyle=*](0,1)(-3.5,0)(-2.5,0)(-0.5,0)(0.5,0)(1.5,0)(3.5,0)
	\end{pspicture}
	\end{center}
\end{cor}

\proof The result of turning a poset upside-down is to reverse all its linear extensions,
which changes ascents to descents and vice-versa.  So if $F(x)$ is the descent polynomial of
the original poset, the descent polynomial of the new poset is $x^{n-1}F(x^{-1})$.  So the roots
of the new polynomial are the inverses of the roots of the original. \eop

\section{General Behavior} \label{sec:unimodal}

We can say in general how the sequence $\euler{n}{d}_n,\euler{n}{d}_{n-1},\ldots,\euler{n}{d}_1$ behaves.


The set of numbers $\euler{n}{d}_k$, for $n$ fixed, is very nearly unimodal if arranged appropriately.

\begin{theorem}\label{thm:unimodal} Fix $n$ and $d$.  Then

\begin{tabular}{rll}
  (i) & If $d=0$,                       & $0 = \euler{n}{d}_n = \cdots = \euler{n}{d}_2 < \euler{n}{d}_1 = 1$ \\
 (ii) & If $1 \le d \le (n-3)/2$,       & $\euler{n}{d}_n < \euler{n}{d}_{n-1} < \cdots < \euler{n}{d}_1$ \\
(iii) & If $n$ is even and $d=(n-2)/2$, & $\euler{n}{d}_n < \cdots < \euler{n}{d}_2 = \euler{n}{d}_1$ \\
 (iv) & If $n$ is odd and $d=(n-1)/2$,  & $\euler{n}{d}_n < \cdots < \euler{n}{d}_{(n+1)/2} > \cdots > \euler{n}{d}_1$ \\
  (v) & If $n$ is even and $d=n/2$,     & $\euler{n}{d}_n = \euler{n}{d}_{n-1} > \cdots > \euler{n}{d}_1$ \\
 (vi) & If $(n+1)/2 \le d \le n-2$,     & $\euler{n}{d}_n > \euler{n}{d}_{n-1} > \cdots > \euler{n}{d}_1$ \\
(vii) & If $d=n-1$,                     & $1 = \euler{n}{d}_n > \euler{n}{d}_{n-1} = \cdots = \euler{n}{d}_1 = 0$. \\
\end{tabular}
\end{theorem}

\proof (i) follows from the fact that the identity is the only permutation with 0 descents.
(v), (vi), and (vii) follow from (iii), (ii), and (i) respectively because
$\euler{n}{d}_k = \euler{n}{n-1-d}_{n+1-k}$.

Let $f_n(x) = \euler{n}{\floor{x/n} + 1}_{n\floor{x/n}+n-x}$, which means 
that $f_n(nd-k) = \euler{n}{d}_k$ if $0 \le d \le n-1$ and $1 \le k \le n$. 
\begin{figure}
	\centering
	\includegraphics[height=4in,width=4in]{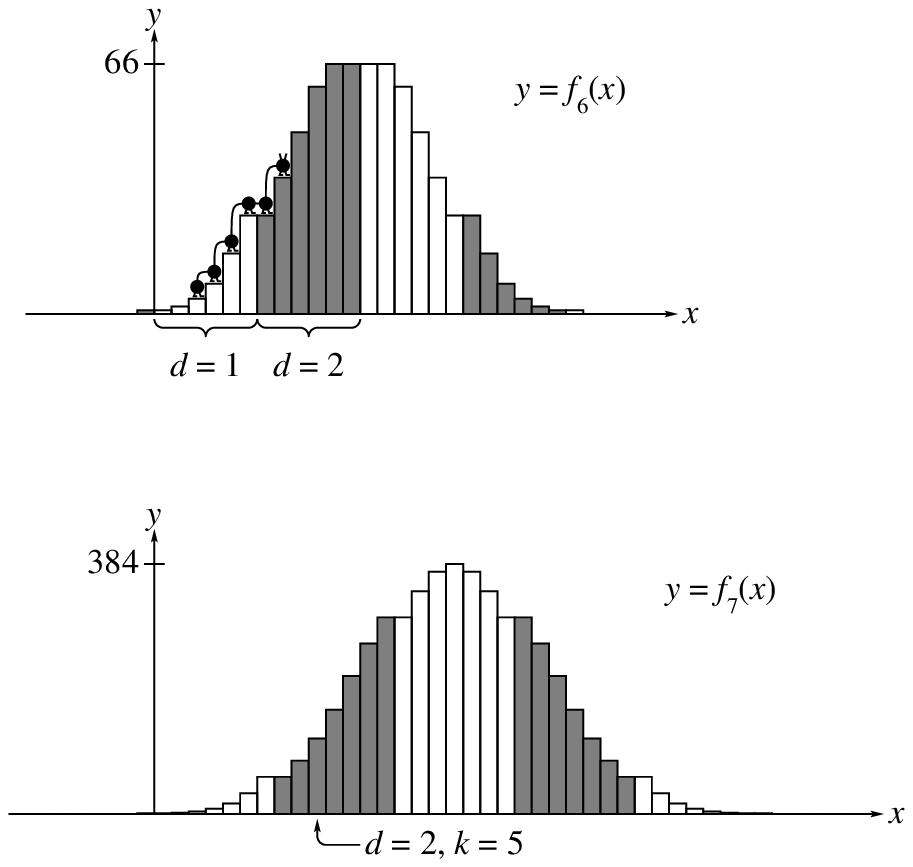}
	\caption{The graphs of $f_6(x)$ and $f_7(x)$, where $f_n(nd-k) = \euler{n}{d}_k$, as defined in \thmnum{unimodal}.}
	\label{fig:graphs}
\end{figure}
Figure~\ref{fig:graphs} shows the graphs of $f_6(x)$ and $f_7(x)$.  Each 
monochromatic section is a sequence of the form $\euler{n}{d}_n, 
\euler{n}{d}_{n-1}, \ldots, \euler{n}{d}_1$.  Note the graphs plateau where 
one sequence meets the next.  Since $\euler{n}{d}_1 = \euler{n-1}{d} = 
\euler{n}{d+1}_n$, each sequence begins where the previous one ends.  The 
content of the theorem is that $f_n$ is basically unimodal.  That is, the 
sequences on the left increase, those on the right decrease, and those in 
the middle behave according to (iii) through (v).

The theorem is true for small $n$ by inspection.  By \eqnum{rec1},
\[ \euler{n+1}{d}_k = \sum_{\ell=1}^{k-1} \euler{n}{d-1}_\ell + \sum_{\ell = k}^n \euler{n}{d}_\ell
   = \sum_{\ell=1}^{k-1} f_n(n(d-1)-\ell) + \sum_{\ell=k}^n f_n(nd-\ell). \]
Let $i = \ell+n-k$ in the first sum and $\ell-k$ in the second and we have
\[ \euler{n+1}{d}_k = \sum_{i=n-k+1}^{n-1} f(n(d-1)-(i-n+k)) + \sum_{i=0}^{n-k} f(nd - (i+k))
   = \sum_{i=0}^{n-1} f(nd - k - i). \]
So imagine a caterpillar of length $n$ crawling on the graph of $y = f_n(x)$, as shown in
the top graph of Figure~\ref{fig:graphs}.  If his head
is at $x$-position $nd-k$, the equation above says that the sum of the heights of his segments
(or his total potential energy) is $\euler{n+1}{d}_k$.  If he were to take a step forward, his total energy
would be $\euler{n+1}{d}_{k+1}$.  That would be an increase in energy if the new height of his head is higher 
than the current height of his tail.  The theorem now follows easily by induction. \eop

\section{Behavior if $d \ll n$} \label{sec:geometric}

If $d$ is much less than $n$, and $\pi$ is selected at random from those permutations of
$n$ letters which have $d$ descents, then the distribution of $\pi(1)$ 
approaches a geometric distribution uniformly, in the following sense.

\begin{theorem} Fix an integer $d > 0$.  Suppose $\pi_n$ is chosen
uniformly from those permutations of $n$ letters which have $d$ descents.  Then
for any $\epsilon > 0$ there is an $N$ such that
\begin{equation}
	\abs{\frac{\pr(\pi_n(1) = k)}{(1-p)p^{k-1}} - 1} < \epsilon
	\label{eq:geomthm}
\end{equation}
for all integers $n$ and $k$ with $n \ge N$ and $1 \le k \le n$, where $p = \frac{d}{d+1}$.
\end{theorem}
\proof For $0 \le j \le d$, let $P_j(n) = (-1)^{d-j}\binom{n}{d-j}$.  Then by \eqnum{andkformula} and \eqnum{eulerformula},
\begin{align}
	\euler{n}{d}_k &= d^{k-1}(d+1)^{n-k} \sum_{0 \le j \le d} P_j(n)\paren{\frac{j}{d}}^{k-1}\paren{\frac{j+1}{d+1}}^{n-k} \\
	\euler{n}{d}   &= (d+1)^n \sum_{0 \le j \le d} P_j(n+1)\paren{\frac{j+1}{d+1}}^n.
\end{align}
Since $(1-p)p^{k-1} = d^{k-1}/(d+1)^k$, the left-hand side of \eqnum{geomthm} is
\[ \abs{\frac{\sum_{0 \le j \le d} P_j(n) \paren{\frac{j}{d}}^{k-1} \paren{\frac{j+1}{d+1}}^{n-k}}
             {\sum_{0 \le j \le d} P_j(n+1) \paren{\frac{j+1}{d+1}}^{n}} - 1} \]
and since $P_d(n) = \binom{n}{0} = 1$, the last term of both sums is 1.  Therefore we have
\[ \abs{\frac{\sum_{0 \le j < d} \bracket{  P_j(n) \paren{\frac{j}{d}}^{k-1} \paren{\frac{j+1}{d+1}}^{n-k}
                                          - P_j(n+1) \paren{\frac{j+1}{d+1}}^n}}
             {1 + \sum_{0 \le j < d} P_j(n+1) \paren{\frac{j+1}{d+1}}^n}}. \]
Since $j/d < (j+1)/(d+1)$ when $0 \le j < d$, that's bounded above by
\[ \frac
	{
		\sum_{0 \le j < d} \bracket{
			\abs{P_j(n)} \paren{\frac{j+1}{d+1}}^{n-1} +
			\abs{P_j(n+1)} \paren{\frac{j+1}{d+1}}^n
		}
	}{
		\abs{1 + \sum_{0 \le j < d} P_j(n+1) \paren{\frac{j+1}{d+1}}^n}
	}.
\]
Now each term in each sum is a polynomial in $n$ times a decaying exponential in $n$.  So both
sums go to 0 as $n$ goes to infinity. \eop
                                    
\begin{cor} The total variation distance between the distribution of $\pi_n(1)$ and the geometric
distribution with parameter $p = \frac{d}{d+1}$ approaches 0 as n approaches infinity. \end{cor}

\section{If Both Ends Are Fixed} \label{sec:bothends}

We might now ask about the number of permutations with $d$ descents whose first and last
positions are fixed.  Let
\[ \euler{n}{d}_k^\ell \assign \sizeof{\set{\pi \in S_n: \mbox{$\des(\pi)=d$, $\pi(1)=k$, and $\pi(n)=\ell$}}}. \]

\begin{theorem} Suppose $1 \le k < k+m \le n$.  Then
\[ \euler{n}{d}_k^{k+m} = \euler{n-1}{d}^m \andd
   \euler{n}{d}_{k+m}^k = \euler{n-1}{d-1}_m. \]
\label{thm:bothends} \end{theorem}

\proof Let $\psi \in S_n$ be the $n$-cycle $(n,n-1,\ldots,2,1)$.  Then for any $\pi \in S_n$,
\[ \psi\pi(i) = \begin{piecewise}
	\pwif{\pi(i)-1}{\pi(i)>1} \\
	\pwif{n}{\pi(i)=1.}
\end{piecewise} \]
(Imagine a device like a car odometer, with a window and $n$ wheels,
on each of which are painted the numbers 1 through $n$.  $\pi$ can be represented
by turning the $i$th wheel until $\pi(i)$ shows through the window, for all $i$.
If one then rolls all the wheels backward a notch, $\psi\pi$ shows
through the window.  For this reason we will refer to the transformation
$\pi \mapsto \psi\pi$ as a {\bfseries rollback}.)

If $1 \le i < n$, let $D_i(\pi) = [\pi(i)>\pi(i+1)]$.  The pair $\pi(i),\pi(i+1)$
has one of four types:

\begin{tabular}{|lr|c|c|c|} \hline
Type & & $D_i(\pi)$ & $D_i(\psi\pi)$ & $D_i(\psi\pi)-D_i(\pi)$ \\ \hline
A & $1<\pi(i)<\pi(i+1)$ & 0 & 0 & 0 \\
B & $1<\pi(i+1)<\pi(i)$ & 1 & 1 & 0 \\
C & $1=\pi(i)<\pi(i+1)$ & 0 & 1 & 1 \\
D & $1=\pi(i+1)<\pi(i)$ & 1 & 0 & -1 \\ \hline
\end{tabular}

Most pairs are of type A or B.  $\pi$ will have one pair of type C unless $\pi(n)=1$
and one pair of type D unless $\pi(1)=1$.  Therefore
\[ \des(\psi\pi)-\des(\pi) = \sum_{i=1}^{n-1} D_i(\psi\pi)-D_i(\pi) = \begin{piecewise}
	\pwif{1}{\pi(1)=1} \\
	\pwif{-1}{\pi(n)=1} \\
	\pwott{0}
\end{piecewise} \]
Let 
\begin{align*}
	P_a^b &\assign \set{\pi\in S_n:\mbox{ $\pi(1)=a$ and $\pi(n)=b$}} \\
	Q^b &\assign \set{\pi\in S_{n-1}: \pi(n)=b}.
\end{align*}
Consider the following sequence of bijections:
\begin{equation*}
\begin{CD}
P_k^{k+m} @>{\mbox{rollback}}>> P_{k-1}^{k-1+m} @>{\mbox{rollback}}>> \cdots @>{\mbox{rollback}}>> P_1^{1+m} @>{\mbox{rollback}}>> P_n^m @>{\mbox{shorten}}>> Q^m
\end{CD}
\end{equation*}
where ``shortening'' a permutation means removing $n$.  (See 
Figure~/ref{fig:rollback} for an example.)  The first $k-1$ rollbacks all 
preserve $\des$, and the final one increments $\des$.  But the shortening 
decrements it again, since it removes $n$ from the front of the permutation.  
Therefore the net effect, across the whole sequence, is to preserve $\des$. 
So $\euler{n}{d}_k^{k+m} = \euler{n-1}{d}^m$ for all $d$.

The second part of the theorem follows from the bijective sequence
\begin{equation*}
\begin{CD}
P_{k+m}^k @>{\mbox{rollback}}>> P_{k-1+m}^{k-1} @>{\mbox{rollback}}>> \cdots @>{\mbox{rollback}}>> P_{1+m}^1 @>{\mbox{rollback}}>> P_m^n @>{\mbox{shorten}}>> Q_m
\end{CD}
\end{equation*}
where $Q_a = \set{\pi\in S_n:\pi(1)=a}$.  Here the final rollback decrements $\des$, and the shortening leaves it
unchanged.  So $\euler{n}{d}_{k+m}^k = \euler{n-1}{d-1}_m$. \eop

\begin{figure} \label{fig:rollback}
	\centering
	\includegraphics[height=170pt,width=470pt]{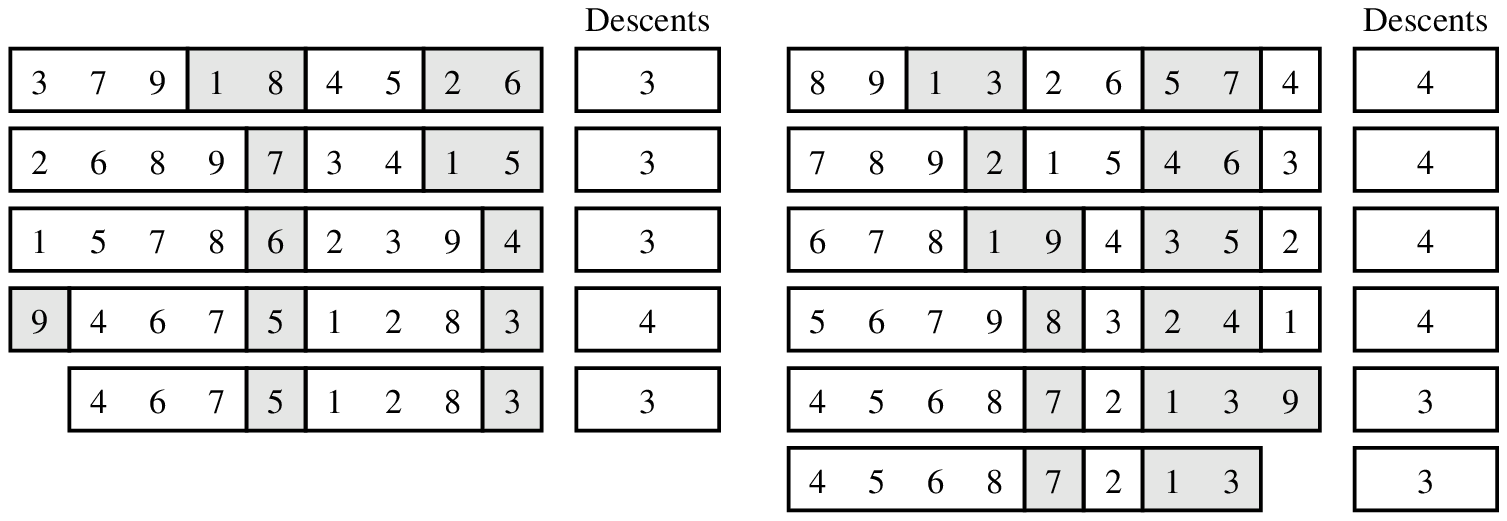}
	\caption{Examples of the actions of the bijections described in Theorem~\ref{thm:bothends}, for $n=9$.  Vertical
	lines show the positions of descents.  If $\pi(1)<\pi(n)$, as at the top left, then
	the permutation is ``rolled back'' until $n$ appears at the front, and then $n$ is removed.  In each of the
	rollbacks but the last, one of the internal bars moves one position to the right, to accomodate a 1 changing to
	an $n$, but the total number of descents stays the same.  Only when the number in the first position changes
	from 1 to $n$ do we gain a descent, but it vanishes again when we remove $n$ in the last step.  The procedure is
	is similar when $\pi(1)>\pi(n)$, as on the right, but the last rollback eliminates a descent, and removing $n$
	leaves the number of descents unchanged.}
\end{figure}

\begin{cor}
	If $1 \le k,\ell \le n$ and $P_{n,k}^\ell$ is the poset on $\set{1,2,\ldots,n}$ defined by
	$k <_P a <_P \ell$ for all $a$ other than $k$ and $\ell$, then the descent polynomial of
	$\mathcal{L}(P_{n,k}^\ell)$ has only distinct real zeroes.
\end{cor}

\proof If $\ell = k+m$, then the polynomial in question is
\[ \sum_d \euler{n}{d}_k^\ell x^d = \sum_d \euler{n-1}{d}^m x^d \]
which was shown to have real distinct zeroes in \cornum{upsidedown}.  As in that
corollary, it follows immediately that turning the poset upside-down inverts the roots
of the polynomial, leaving them real. \eop


\section{Remarks} \label{sec:remarks}

In \secnum{stein} we noted that if $\pi$ is uniformly distributed over $S_n$ then
the distribution of $D = \des(\pi)$ is approximately normal.  Thus the normal density function
\[ \frac1{\sqrt{2\pi}\sigma} \exp\curly{-\frac12\paren{\frac{d-\mu}{\sigma}}^2} \]
with $\mu = \frac{n-1}2$ and $\sigma = \sqrt{\frac{n+1}{12}}$ is a good approximation for $\frac1{n!} \euler{n}{d}$
when $d$ is close to $\mu$.  However, it can be off by orders of magnitude when $d$ is very small or very large.

Theorem~\ref{thm:unimodal} shows that the sequence $\euler{n}{d}_1, \euler{n}{d}_2, \ldots, \euler{n}{d}_n$ 
is an interpolation between $\euler{n-1}{d}$ and $\euler{n-1}{d-1}$, so it seems a reasonable hypothesis that
if $d$ is close to $\frac{n-1}2$, then $\euler{n}{d}_k$ is well approximated by
\[ \frac{(n-1)!}{\sqrt{\frac{\pi n}6}} \exp\curly{-\frac12 \paren{\frac{d + \frac{n-k}{n-1} - \frac{n}2}{\frac{\sqrt{n}}{12}}}^2}. \]
Experimental evidence for $n\le 200$ suggests that this is in fact the case.  So while the distribution
of $\pi(1)$ given $\des(\pi)$ is by no means normal, it does seem to behave like a segment of the normal curve
when $d$ is near $\frac{n-1}2$.

More generally, there may be some underlying curve which the 
Eulerian numbers, properly normalized, can be said to approach as $n$ grows 
large.  It will look like a bell curve, but not be exactly normal, since the normal
approximation is not very good when $d\ll n$.  If so, it seems likely that the refined
Eulerian numbers presented in this paper can be said to approach points on the same curve.

\section{Acknowledgements}

Sergi Elizalde found an alternate version of the generating function in 
\thmnum{gfall}, and also sketched an alternate proof for \thmnum{edpi1} 
using generating functions.  The author would also like to thank Persi 
Diaconis and Jason Fulman for introductions to Stein's method, and Sergey 
Fomin, Mark Skandera, and Boris Pittel for explanations of the Neggers-Stanley
problem.  Bruce Sagan and Herbert Wilf encouraged this work forward. 
And the most thanks go to Divakar Viswanath, who discussed every aspect of 
this paper at every stage of developement.

\nocite{ConcreteMath}
\nocite{Wilf}

\bibliography{ref}

\end{document}

%% file: defs.tex
\newenvironment{piecewise}{\left\{\begin{array}{ll}}{\end{array}\right.}
\newcommand{\pwif}[2]{#1 & \mbox{if \ensuremath{#2}}}
\newcommand{\pwot}[1]{#1 & \mbox{otherwise}}
\newcommand{\pwott}[1]{#1 & \mbox{otherwise.}}

\newcommand{\ex}{\mathbb{E}}
\newcommand{\pr}{\mathbb{P}}
\newcommand{\var}{\mathop{\mathrm{Var}}}

\newcommand{\abs}[1]{\ensuremath{\left| #1 \right|}}
\newcommand{\sizeof}[1]{\##1}
\newcommand{\assign}{\ensuremath{:=}}
\newcommand{\del}[2]{\ensuremath{\frac{\partial#1}{\partial#2}}}
\newcommand{\set}[1]{\left\{#1\right\}}
\newcommand{\paren}[1]{\ensuremath{\left( #1 \right)}}
\newcommand{\bracket}[1]{\ensuremath{\left[ #1 \right]}}
\newcommand{\curly}[1]{\ensuremath{\left\{ #1 \right\}}}
\newcommand{\floor}[1]{\ensuremath{\left\lfloor#1\right\rfloor}}
\newcommand{\ceil}[1]{\ensuremath{\left\lceil#1\right\rceil}}
\newcommand{\andd}{\mbox{\quad and\quad}}
\newcommand{\orr}{\mbox{\quad or\quad}}
\newcommand{\smallhead}[1]{{\bfseries\textsf{#1}: }}
\newcommand{\definition}{\smallhead{Definition}}
\newcommand{\notation}{\smallhead{Notation}}
\newcommand{\note}{\smallhead{Note}}
\newcommand{\ha}[1]{{\Huge\bfseries\textsf{#1}}}
\newcommand{\hb}[1]{{\huge\bfseries\textsf{#1}}}
\newcommand{\hc}[1]{{\LARGE\bfseries\textsf{#1}}}
\newcommand{\hd}[1]{{\Large\bfseries\textsf{#1}}}
\newcommand{\he}[1]{{\large\bfseries\textsf{#1}}}

\newcommand{\euler}[2]{\genfrac{\langle}{\rangle}{0pt}{}{#1}{#2}}

\newcommand{\des}{\mathop{\mathrm{Des}}}
\newenvironment{myitemize}{\begin{list}{$\bullet$}{\setlength{\topsep}{0pt}\setlength{\parsep}{0pt}}}{\end{list}}
\newcommand{\evalat}[2]{\ensuremath{\left.#1\right|_{#2}}}
\newcommand{\eqnum}[1]{Equation~(\ref{eq:#1})}
\newcommand{\thmnum}[1]{Theorem~\ref{thm:#1}}
\newcommand{\lemnum}[1]{Lemma~\ref{lem:#1}}
\newcommand{\cornum}[1]{Corollary~\ref{cor:#1}}
\newcommand{\secnum}[1]{Section~\ref{sec:#1}}

\newtheorem{theorem}{Theorem}
\newtheorem{lemma}[theorem]{Lemma}
\newtheorem{cor}[theorem]{Corollary}
\newcommand{\proof}{\textit{Proof.}\hspace{1em}} 
\newcommand{\eop}{\hfill\ensuremath{\Box}} 

\newenvironment{prototheorem}[1]{
	\renewcommand{\thetheorem}{\ref{thm:#1}}
	\begin{theorem}
}{
	\end{theorem}
	\renewcommand{\thetheorem}{\arabic{theorem}}
	\addtocounter{theorem}{-1}
}